\documentstyle[11pt]{article}
 \input{psfig}

\newcommand{\qed}{\mbox{$\Diamond$}\vspace{\baselineskip}}

\newtheorem{theorem}{Theorem}

\newtheorem{lemma}{Lemma}
\newtheorem{proposition}{Proposition}
\newtheorem{definition}{Definition}

\newtheorem{example}{Example}
 
\newenvironment{proof}{\noindent {\bf Proof:}}{{\qed}}

\newcommand{\vanish}[1]{}
 
\begin{document}
\title{A self-dual poset on objects counted by the Catalan numbers}
 
\author{Mikl\'os B\'ona\thanks{This paper was written while the author's stay at the
Institute was supported by Trustee Ladislaus von Hoffmann, the Arcana
Foundation.}\\
        School of Mathematics\\
        Institute for Advanced Study\\
        Princeton, NJ 08540 \\
        }
\maketitle
\begin{abstract}
We examine the poset $P$ of 132-avoiding $n$-permutations ordered by
descents. 
We show that this poset is the "coarsening" of the well-studied poset $Q$ of
 noncrossing partitions . In other words, if
$x<y$ in $Q$, then $f(y)<f(x)$ in $P$, where $f$ is 
the canonical bijection from the set of noncrossing partitions
onto that of 
132-avoiding permutations. This enables us to prove many properties of $P$.
\end{abstract}

\section{Introduction}
There are more than 150 different objects enumerated by Catalan numbers. Two
of the most carefully studied ones are noncrossing partitions  and 
132-avoiding permutations. A partition $\pi=(\pi_1,\pi_2,\cdots ,\pi_t)$
  of the set $[n]=\{1,2,\cdots, n\}$ is called noncrossing \cite{kreweras}
 if it has no four 
elements $a<b<c<d$ so 
that $a,c\in \pi_i$ and $b,d \in \pi_j$
for some distinct $i$ and $j$. A permutation of $[n]$, or, in what follows,
an $n$-permutation, is called 132-avoiding
\cite{schmidt} if it does not have three entries $a<b<c$ so that $a$ is
 the leftmost of them and $b$ is the rightmost of them.

Noncrossing partitions of $[n]$
 have a natural and well studied partial order: the
refinement order
$Q_n$. In this order $\pi_1<\pi_2$ if each block of $\pi_2$ is the union of
some blocks of $\pi_1$. The poset $Q_n$ is known to be a lattice,
 and it is graded, rank-symmetric, rank-unimodal, and $k$-Sperner 
\cite{edelman}. The poset $Q_n$ has been proved to be self-dual in two steps
  \cite{kreweras}, \cite{simull}.  

In this paper we introduce a new partial order of 132-avoiding $n$-permutations
which will naturally translate into one of noncrossing partitions. In this
poset, for two 132-avoiding $n$-permutations $x$ and $y$,
 we define $x<y$ if the
descent set of $x$ is contained in that of $y$. (We will provide a natural
 equivalent, definition, too.) 
We will see that this new partial order $P_n$ is a coarsening
 of the dual of $Q_n$.
In other words, if for two noncrossing partitions $\pi_1$ and $\pi_2$ we have 
$\pi_1<\pi_2$ in $Q_n$, then we also have $f(\pi_2)<f(\pi_1)$ in $P_n$,
 where $f$ 
is a natural bijection from the set of noncrossing partitions onto that
of 132-avoiding permutations. This will enable us to prove that $P_n$ has the
same rank-generating function as $Q_n$, and so
$P_n$ is rank-unimodal, rank-symmetric and $k$-Sperner. Furthermore, we will
also prove that $P_n$ is self-dual in a somewhat more direct way than it is
proved for $Q_n$.

\section{Our main results}
\subsection{A bijection and its properties}
It is not difficult to find a bijection from the set of noncrossing partitions
of $[n]$ onto that of 132-avoiding $n$-permutations. However, we will
exhibit such a bijection here and analyze its structure as it will be
our major tool in proving our theorems. To avoid confusion, integers
belonging to a partition will be called {\em elements}, while integers 
belonging to a permutation will be called {\em entries}. An $n$-permutation
$x=x_1x_2\cdots x_n$ will always be written in the one-line notation, with
$x_i$ denoting its $i$th entry.

Let $\pi$ be a noncrossing partition of $[n]$. We construct the 132-avoiding
permutation $p=f(\pi)$ corresponding to it. Let $k$ be the largest element of
 $\pi$ which is in
the same block of $\pi$ as 1. Put the entry $n$ of $p$ to the $k$th 
position, so $p_k=n$. As $p$ is to be 132-avoiding, 
this implies that entries larger
than $n-k$ are on the left of $n$ and entries less than or equal to $n-k$ 
are on the right of $n$ in $q$.

Then we continue this procedure recursively. As $\pi$ is noncrossing, blocks
which contain elements larger than $k$ cannot contain elements smaller than 
$k$. Therefore, the restriction of $\pi$ to $\{k+1,k+2,\cdots ,n\}$ is a 
noncrossing partition, and it corresponds to the 132-avoiding permutation
of $\{1,2,\cdots n-k\}$ which is on the left of $n$ in $\pi$ by this same
 recursive procedure.

We still need to say what to do with blocks of $\pi$ containing elements
smaller than or equal to $k$. Delete $k$, and apply this same procedure
for the resulting noncrossing partition on $k-1$ elements. This way we
obtain a 132-avoiding permutation of $k-1$ elements, and this is what we
needed for the part of $p$ on the left of $n$, that is, for $\{n-k+1,
n-k+2,\cdots ,n-1\}$.

So in other words, if $\pi_1$ is the restriction of $\pi$ into $[k-1]$ and
$\pi_2$ is the restriction of $\pi$ into $\{k+1,k+2,\cdots ,n\}$, then $f(\pi)$
is the concatenation of $f(\pi_1)$, $n$ and $f(\pi_2)$, where $f(\pi_1)$
 permutes
the set $\{n-k+1,n-k+2,\cdots ,n-1\}$  and $f(\pi_2)$ permutes the set $[n-k]$.

To see that this is a bijection note that we can recover the largest
element of the block containing the entry  1 from the position of $n$ in $p$
and then proceed recursively. 

\begin{example} {\em If $\pi=(\{1,4,6\}, \{2,3\}, \{5\}, \{7,8\})$, then 
$f(\pi)=64573812$. }\end{example}

\begin{example} {\em If $p=(\{1,2,\cdots ,n\})$, then 
$f(p)=12\cdots n $. }\end{example}

\begin{example} {\em If $p=(\{1\},\{2\},\cdots ,\{n\})$, then 
$f(p)=n\cdots 21 $. }\end{example}

The following definition is widely used in the literature.
\begin{definition} Let $p=p_1p_2\cdots p_n$ be a permutation. We say the $i$
is a {\em descent} of $p$ if $p_i>p_{i+1}$. The set of all descents of $p$
is called the {\em descent set} of $p$ and is denoted $D(p)$. \end{definition}

Now we are in a position to define the poset $P_n$ of 132-avoiding
permutations we want to study. 
\begin{definition} 
Let $x$ and $y$ be two 132-avoiding $n$-permutations. 
We say that $x<_Py$ (or $x<y$
in $P_n$) if 
$D(x)\subset D(y)$. \end{definition}
Clearly, $P_n$ is a poset as inclusion is transitive. It is easy to see that in
132-avoiding permutations, $i\leq 1$
 is a descent if and only if $p_{i+1}$ is smaller
than every entry on its left, (such an element is called a left-to-right
minimum). So $x<_Py$ if and only if the set of positions in which $x$ has
a left-to-right minimum is a proper subset of that of those positions in
which $y$ has a left-to-right minimum.
The Hasse diagram of $P_4$ is shown on the Figure below.
\vskip .3 cm
\begin{center}
\setlength{\unitlength}{0.00033300in}%
\begingroup\makeatletter\ifx\SetFigFont\undefined
\def\x#1#2#3#4#5#6#7\relax{\def\x{#1#2#3#4#5#6}}%
\expandafter\x\fmtname xxxxxx\relax \def\y{splain}%
\ifx\x\y   
\gdef\SetFigFont#1#2#3{%
  \ifnum #1<17\tiny\else \ifnum #1<20\small\else
  \ifnum #1<24\normalsize\else \ifnum #1<29\large\else
  \ifnum #1<34\Large\else \ifnum #1<41\LARGE\else
     \huge\fi\fi\fi\fi\fi\fi
  \csname #3\endcsname}%
\else
\gdef\SetFigFont#1#2#3{\begingroup
  \count@#1\relax \ifnum 25<\count@\count@25\fi
  \def\x{\endgroup\@setsize\SetFigFont{#2pt}}%
  \expandafter\x
    \csname \romannumeral\the\count@ pt\expandafter\endcsname
    \csname @\romannumeral\the\count@ pt\endcsname
  \csname #3\endcsname}%
\fi
\fi\endgroup
\begin{picture}(11333,8460)(301,-7828)
\thicklines
\put(901,-2536){\circle{300}}
\put(2701,-2536){\circle{300}}
\put(4276,-2536){\circle{300}}
\put(6751,-2461){\circle{300}}
\put(9001,-2461){\circle{300}}
\put(11476,-2461){\circle{300}}
\put(4501,-5086){\circle{300}}
\put(6707,-4980){\circle{300}}
\put(9032,-5055){\circle{300}}
\put(11251,-5086){\circle{300}}
\put(5776,-7336){\circle{300}}
\put(5326,314){\circle{300}}
\put(2820,-5117){\circle{300}}
\put(826,-5161){\circle{300}}
\put(5251,314){\line(-3,-2){4378.846}}
\put(5251,314){\line(-5,-6){2446.721}}
\put(5251,314){\line(-1,-3){975}}
\put(5251,314){\line( 1,-2){1440}}
\put(5251,314){\line( 4,-3){3732}}
\put(5326,314){\line( 2,-1){6030}}
\put(826,-5236){\line( 5,-2){5017.241}}
\put(5776,-7411){\line( 5, 2){5469.828}}
\put(5776,-7411){\line( 4, 3){3132}}
\put(2776,-5086){\line( 5,-4){2963.415}}
\put(4501,-5011){\line( 1,-2){1215}}
\put(6751,-4936){\line(-2,-5){987.931}}
\put(826,-2536){\line( 0,-1){2700}}
\put(826,-2536){\line( 4,-5){2034.146}}
\put(826,-2611){\line( 3,-2){3738.461}}
\put(826,-2611){\line( 5,-2){5909.483}}
\put(2701,-2611){\line(-3,-4){1935}}
\put(2701,-2611){\line( 0,-1){2550}}
\put(2851,-5161){\line( 0, 1){ 75}}
\put(2701,-2611){\line( 3,-4){1863}}
\put(2776,-2611){\line( 5,-3){3948.529}}
\put(2776,-2611){\line( 5,-2){6193.966}}
\put(4276,-2611){\line( 1,-6){407.432}}
\put(4276,-2686){\line( 6,-5){2567.213}}
\put(4276,-2686){\line( 2,-1){4710}}
\put(901,-5161){\line( 5, 2){5948.276}}
\put(2851,-5086){\line( 3, 2){3876.923}}
\put(6751,-2536){\line(-5,-6){2176.229}}
\put(9076,-2611){\line(-2,-1){4560}}
\put(9076,-2611){\line(-5,-2){6219.828}}
\put(9076,-2611){\makebox(16.6667,25.0000){\SetFigFont{10}{12}{rm}.}}
\put(11476,-2536){\line(-2,-1){4740}}
\put(11401,-2611){\line(-1,-1){2400}}
\put(9001,-2461){\line(-3,-1){8167.500}}
\put(816,-5148){\line( 4, 3){3468}}
\put(899,-2617){\line( 3,-1){7942.500}}
\put(8966,-2458){\line( 5,-6){2145.492}}
\put(11141,-5008){\line( 0, 1){  0}}
\put(11141,-5008){\line( 1, 6){395.270}}
\put(6668,-2401){\line( 5,-3){4577.206}}
\put(4201,-2611){\line(-1,-2){1260}}
\put(2851,-5086){\line( 1, 2){1260}}
\put(4951,464){\makebox(0,0)[lb]{\smash{\SetFigFont{8}{9.6}{rm}4321}}}
\put(5476,-7786){\makebox(0,0)[lb]{\smash{\SetFigFont{8}{9.6}{rm}1234}}}
\put(301,-2311){\makebox(0,0)[lb]{\smash{\SetFigFont{8}{9.6}{rm}3214}}}
\put(2101,-2311){\makebox(0,0)[lb]{\smash{\SetFigFont{8}{9.6}{rm}4213}}}
\put(3676,-2311){\makebox(0,0)[lb]{\smash{\SetFigFont{8}{9.6}{rm}4312}}}
\put(5926,-2311){\makebox(0,0)[lb]{\smash{\SetFigFont{8}{9.6}{rm}4231}}}
\put(8851,-2236){\makebox(0,0)[lb]{\smash{\SetFigFont{8}{9.6}{rm}3241}}}
\put(11176,-2236){\makebox(0,0)[lb]{\smash{\SetFigFont{8}{9.6}{rm}3421}}}
\put(2401,-5536){\makebox(0,0)[lb]{\smash{\SetFigFont{8}{9.6}{rm}3124}}}
\put(376,-5536){\makebox(0,0)[lb]{\smash{\SetFigFont{8}{9.6}{rm}2134}}}
\put(6001,-5311){\makebox(0,0)[lb]{\smash{\SetFigFont{8}{9.6}{rm}3412}}}
\put(3976,-5461){\makebox(0,0)[lb]{\smash{\SetFigFont{8}{9.6}{rm}4123}}}
\put(10501,-5236){\makebox(0,0)[lb]{\smash{\SetFigFont{8}{9.6}{rm}2341}}}
\put(8026,-5311){\makebox(0,0)[lb]{\smash{\SetFigFont{8}{9.6}{rm}2314}}}
\end{picture}

\end{center}
\vskip .3 cm
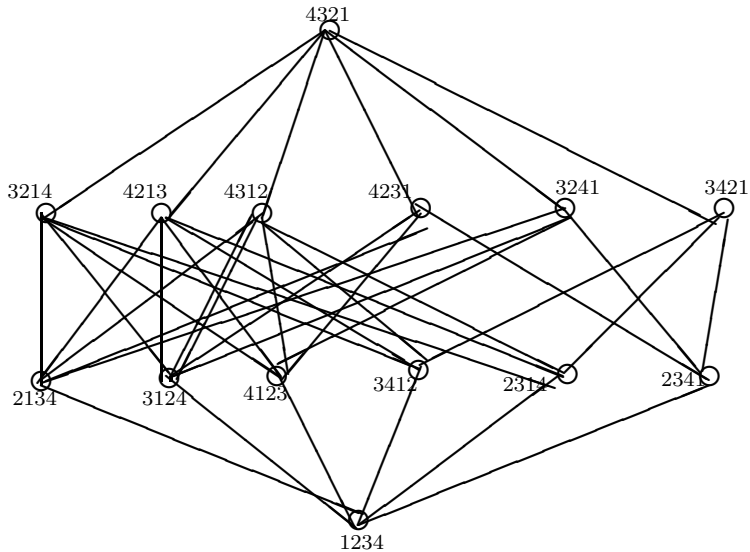
\centerline{Figure 1: The Hasse diagram of $P_4$. }

The following proposition describes the relation between the blocks of
$x$ and the descent set of $f(x)$. 

\begin{proposition} \label{block} The bijection $f$ has the following property:
$i\in D(f(x))$ if and only if $i+1$ is the smallest element of its block.
\end{proposition}
\begin{proof} 
By induction on $n$. For $n=1$ and $n=2$ the statement is true. Now suppose
we know the statement for all positive integers smaller than $n$. 
Then we distinguish two cases:
\begin{enumerate} \item If 1 and $n$ are in the same block of $x$, then the 
construction of $f(x)$ simply starts by putting the entry $n$ to the
last slot of $f(x)$, then deleting the element  $n$ from $x$. Neither of these
steps alters the set of minimal elements of blocks or that of 
descents in any way. Therefore, the algorithm is reduced to one of size
$n-1$, and the proof follows by induction.
 \item If the largest element $k$ of the block containing 1 is smaller than
$n$, then as we have seen above,  $f$ constructs the images of $x_1$ and $x_2$
which will be separated by the entry $n$. Therefore, by the induction 
hypothesis, the descents of $f(x)$ are given by the minimal elements of the
blocks of $x_1$ and $x_2$, and these are exactly the blocks of $x$. There
will also be a descent at $k$ (as the entry $n$ goes to the $k$th slot), and
 that is in accordance with our statement as $k+1$ is certainly the smallest
element of its block. 
\end{enumerate}
\end{proof}

We point out that this implies that $P_n$ is equivalent to a poset of 
noncrossing partitions in which $\pi_1<\pi_2$ if the set of elements which 
are minimal in 
their block in $\pi_1$ is contained in that of elements which are minimal in
their block in $\pi_2$.
\subsection{Properties of $P_n$}
Now we can prove the main result of this paper.
\begin{theorem}
The poset $P_n$ is coarser than the dual of the poset $Q_n$ of noncrossing 
 partitions ordered by refinement. That is, if $x<y$ in $Q_n$, then $f(y)<f(x)$
 in $P_n$.
\end{theorem}
\begin{proof}
If $x<y$, then each block of $x$ is a subset of a block of $y$. Therefore,
if $z$ is the minimal element of a block $B$ of $y$, then it is also the
 minimal
element of the block $E$ of $x$ containing it as  $E\subseteq B$. Therefore,
the set of elements which are minimal in their respective blocks in $x$ 
contains that of
elements which are minimal in their respective blocks in $y$. By Proposition
\ref{block} this implies $D(f(y))\subset D(f(x))$.
\end{proof}

Now we apply this result to prove some properties of $P_n$. For definitions,
see \cite{stanley}.
\begin{theorem} The rank generating function of $P_n$ is equal to that
 of $Q_n$.
 In particular, $P_n$ is rank-symmetric, rank-unimodal and $k$-Sperner.
\end{theorem}
\begin{proof}
By proposition \ref{block}, the number of 132-avoiding permutations having
 $k$ descents equals that of
noncrossing partitions having $k$ blocks, and this is known to be the $(n,k)$ 
Narayana-number $\frac{1}{n}\cdot {n\choose k}{n\choose k-1}$. Therefore
$P_n$ is graded, 
rank-symmetric and rank-unimodal, and its rank generating function is the same
as that of $Q_n$, as $Q_n$ too is graded by the number of blocks 
(and is self-dual).
 As $P_n$ is coarser than $Q_n$, any 
antichain of $P_n$ is an antichain of $Q_n$,
 and the $k$-Sperner property follows.
\end{proof}

We need more analysis to prove that $P_n$ is self-dual, that is, that $P_n$ is
invariant to ``being turned upside down''. Denote $Perm_n(S)$ the
number of 132-avoiding $n$-permutations with descent set $S$. The following
lemma is the base of our proof of self-duality. For $S\subseteq [n-1]$, we
define  $\alpha(S)$ to be the ``reverse complement'' of $S$, that is,
 $i\in \alpha(S)
\Longleftrightarrow n-i \notin S$.
\begin{lemma} For any $S\subseteq [n-1]$, we have
 $Perm_n(S)=Perm_n(\alpha(S))$.
\end{lemma} 
\begin{proof}By induction on $n$. For $n=1,2,3$ the statement is true.
Now suppose we know it for all positive integers smaller than $n$. 
Denote $t$ the smallest element of $S$.
\begin{enumerate} \item Suppose that $t>1$. 
This means
that $x_1<x_2<\cdots <x_t$, and that $x_1,x_2,\cdots ,x_t$ are {\em consecutive
integers}. Indeed, if there were a gap among them, that is, there were an 
integer $y$ so that $y\neq x_i$ for $1\leq i \leq t$, while $x_1<y<x_t$, then
$x_1\:x_t\:y$ would be a 132-pattern. So once we know $x_1$,
 we have noly one choice for $x_2,x_3,\cdots ,x_t$. This implies 
\begin{equation} \label{1a}
Perm_n(S)=Perm_{n-(t-1)}(S-(t-1)),\end{equation} 
 where $S-(t-1)$ is the set obtained from
$S$ by subtracting $t-1$ from each of its elements. 

On the other hand, we have $n-t+1,n-t+2,\cdots ,n-1 \in 
\alpha(S)$, meaning that
$x_{n-t+1}>x_{n-t+2}>\cdots >x_n$, and also, we must have $(
x_{n-t+2},\cdots ,x_n)=(t-1,t-2,\cdots 1)$, otherwise a 132-pattern is formed.
Therefore, \begin{equation}\label{1b}
Perm_n(\alpha(S))=Perm_{n-(t-1)}(\alpha(S)|n-(t-1))\end{equation}
where $\alpha(S)|n-(t-1)$ is simply $\alpha(S)$ without its last $t-1$ 
elements. Clearly, $Perm_{n-(t-1)}(S-(t-1))=Perm_{n-(t-1)}(\alpha(S)|n-(t-1))$
by the induction hypothesis, so equations (\ref{1a}) and  (\ref{1b}) imply
 $Perm_n(S)=Perm_n(\alpha(S))$.
\item If $t=1$, but $S\neq [n-1]$,
 then let $u$ be the smallest index which is not
 in $S$. Then again, $x_u$ must be the smallest positive integer $a$ which
is larger than $x_{u-1}$ and is not equal to some $x_i$, $i\leq u-1$, 
otherwise $x_{u-1}\:x_u\:a$ would be a 132-pattern. So again, we have only one
choice for $x_u$.
On the other hand, the largest index in $\alpha(S)$ will be $n-(u-1)$.  Then
as above, we will only have once choice for $x_{n-u}$. Now we can delete
$u$ from $S$ and $n-u$ from $\alpha(S)$ and proceed by the induction
hypothesis as in the previous case.
\item Finally, if $S=[n-1]$, then the statement is trivially true as 
$Perm_n(S)=Perm_n(\alpha(S))=1.$
\end{enumerate}
So we have seen that $Perm_n(S)=Perm_n(\alpha(S))$ in all cases.
\end{proof}

Now we are in position to prove our next theorem.
\begin{theorem} The poset $P_n$ is self-dual.
\end{theorem}
\begin{proof}
It is clear that in $P_n$  permutations with the same descent set will
cover the same elements and they will be covered by the same elements.
Therefore, such permutations form orbits of $Aut(P_n)$ and they can be permuted
among each other arbitrarily by elements of $Aut(P_n)$. 
One can think of $P_n$ as a Boolean algebra $B_{n-1}$ in which 
some elements have
several copies. One natural anti-automorphism of a Boolean-algebra is 
``reverse complement'', that is, for $S\subseteq [n-1]$, $i\in \alpha(S)
\Longleftrightarrow n-i \notin S$.
 To show that $P_n$ is
self-dual, it is therefore sufficient to show that the corresponding elements
appear with the same multiplicities in $P_n$. So in other words we must show 
that there are as many 132-avoiding permutations with descent set
$S$ as there are with descent set $\alpha(S)$. And that has been proved in the
Lemma.
\end{proof}
\subsection{Further directions}
It is natural to ask for what related combinatorial objects  we could define
such a natural partial order which would turn out to be self-dual and possibly,
have some other nice properties.  {\em Two-stack
sortable permutations} \cite{doron} are an obvious candidate. 
It is known \cite{schaeffer} 
that there are as many of them with $k$ descents as with
$n-1-k$ descents, however, the poset obtained by the descent ordering is
not self-dual, even for $n=4$, so another ordering is needed. Another candidate
could be the poset of the 
recently introduced noncrossing partitions for classical
reflection groups \cite{reiner}, some of which are self-dual in the traditional
refinement order. 

\end{document}